\documentclass[times,preprint,sort&compress]{elsarticle}
\usepackage{epstopdf}
\usepackage[caption=false]{subfig}
\usepackage[hidelinks]{hyperref}
\usepackage[active]{srcltx}
\usepackage{bbm}
\usepackage{amsthm}
\usepackage{amssymb}
\usepackage{stackrel}
\usepackage{xcolor}
\usepackage{amsmath,amscd}

\newtheorem{dl}{Theorem}
\newtheorem{tl}[dl]{Corollary}
\newtheorem{yl}[dl]{Lemma}
\newtheorem{dy}[dl]{Definition}
\newtheorem{lz}[dl]{Example}

\newtheorem{prop}[dl]{Proposition}

\newtheorem{remark}[dl]{Remark}

\newcommand{\be}{\begin{equation}}
\newcommand{\ee}{\end{equation}}
\newcommand{\ba}{\begin{array}}
\newcommand{\ea}{\end{array}}
\newcommand{\bmn}{\begin{eqnarray}}
\newcommand{\emn}{\end{eqnarray}}
\newcommand{\bnm}{\begin{eqnarray*}}
\newcommand{\enm}{\end{eqnarray*}}
\newcommand{\bln}{\begin{subequations}}
\newcommand{\eln}{\end{subequations}}

\newproof{pot444}{\bf Proof of Theorem \ref{mainthm}}
\newproof{pot555}{\bf Proof of Theorem \ref{dlsecond}}
\newproof{pot666}{\bf Proof of Theorem \ref{dlthird}}
\newproof{pot777}{\bf Proof of Corollary \ref{dlthird-one}}
\newcommand{\poq}[2]{(#1;q)_{#2}}

\def\qed{\hfill \rule{4pt}{7pt}}
\def\pf{\noindent {\it Proof.} }
\newcommand{\poqq}[2]{(#1;q^2)_{#2}}

\begin{document}
\title{Two new Bailey pairs and their $q$-identities of Rogers-Ramanujan type modulo 15, 24, and 30}
\author{Jianan Xu\fnref{fn1}}
\fntext[fn1]{E-mail address: 20224007007@stu.suda.edu.cn}
\address[P.R.China]{Department of Mathematics, Soochow University, SuZhou 215006, P.R.China}
\author{Xinrong Ma\fnref{fn2,fn3}}
\fntext[fn2]{This  work was supported by NSFC grant No. 12471315}
\fntext[fn3]{Corresponding author. E-mail address: xrma@suda.edu.cn.}
\address[P.R.China]{Department of Mathematics, Soochow University, SuZhou 215006, P.R.China}
\markboth{J.N. Xu and  X. R. Ma}{Two new Bailey pairs and their $q$-identities of Rogers-Ramanujan type modulo integers}
\begin{abstract}In this paper, we first establish   two    new  Bailey pairs via finding two generalizations  of  Euler's pentagonal number theorem. Next, we specificize the Bailey lemmas with these two Bailey pairs. As applications, we finally establish some $q$-series transformations and $q$-identities of Rogers-Ramanujan type modulo  15, 24, and 30.\end{abstract}
\maketitle
\parskip 7pt

\section{Introduction}
 It is well known that  all sum-product $q$-identities of the form
\begin{align*}
\sum_{n=0}^{\infty} \frac{q^{n^2}}{(q ; q)_n}=\frac{1}{\left(q, q^4 ; q^5\right)_{\infty}}, \quad \sum_{n=0}^{\infty} \frac{q^{n(n+1)}}{(q ; q)_n}=\frac{1}{\left(q^2, q^3 ; q^5\right)_{\infty}}
\end{align*}
are called to be $q$-identities of  Rogers-Ramanujan type. As discovered by   G. E. Andrews  \cite{andrews2}, these two modular identities can be reproduced via the method of Bailey pairs and Bailey chains. 
The purpose of the present paper is to introduce and study two  peculiar but new Bailey pairs, which may serve as two generalizations of Euler's pentagonal number theorem (cf. \cite{andrewseuler}). More importantly, some $q$-identities of Rogers-Ramanujan type are derived from the corresponding Bailey lemmas. As a prelude to our discussion, we introduce some notations and  basic background material. Unless otherwise stated, we will  adopt the standard notation and terminology for  $q$-series from the book  \cite{10} by G. Gasper and M. Rahman.  As usual, the $q$-shifted factorials of the variable $x$ and the base $q$, $0<|q|<1$, are  given by
\begin{align*}
(x;q)_\infty
=\prod_{n=0}^{\infty}(1-xq^n),~~
(x;q)_n=\frac{(x;q)_\infty}{(xq^n;q)_\infty}
\end{align*}
for all $n\in \mathbb{Z}$, $\mathbb{Z}$ is the set of integers. We also employ
 the multi-parameter compact notation
\[(a_1,a_2,\ldots,a_n;q)_\infty=(a_1;q)_\infty (a_2;q)_\infty\ldots (a_n;q)_\infty.\]
In particular, we define the $q$-binomial coefficient 
\begin{align*}
\genfrac{[}{]}{0pt}{}{n}{k}_{q}:=\left\{\begin{array}{ll}\displaystyle
 \frac{\poq{q}{n}}{\poq{q}{k}\poq{q}{n-k}},&\mbox{if}~~n\geq k\geq 0;\\
 \displaystyle\\
 0, &\mbox{otherwise}.\end{array}\right.
\end{align*}
where $n\geq k$ are nonnegative integers.
Hereafter, for convenience, we introduce the notation $\tau(n)$ for $(-1)^{n}q^{\binom{n}{2}}$.  

Our argument  frequently involves  the classical Jacobi triple product identity.
\begin{yl}[\mbox{Jacobi's triple product identity:\cite[(II.28)]{10}}]
\begin{align}
\sum_{n=-\infty}^\infty\,\tau(n)x^n=\poq{x,q/x,q}{\infty}.\label{triple}
 \end{align}
\end{yl}

We also require the concept of Bailey pairs \cite{andrews2,bailey}.
\begin{dy}[The Bailey pair]\label{defnew-ooo}
A Bailey pair relative to $a$ is conveniently  defined to be a pair of sequences $(\alpha_n(a,q))_{n\geq 0}$ and $(\beta_n(a,q))_{n\geq 0},$  for brevity denoted by  $(\alpha_n(a,q),\beta_n(a,q))$,  satisfying
\begin{subequations}\label{baileypairs}
\begin{align}
\beta_n(a,q)&=\sum^{n}_{k=0} \frac{\alpha_k(a,q)}
{\poq{q}{n-k}\poq{aq}{n+k}}.\label{defnew}
\end{align}
Or, by inverting
\begin{align}
\alpha_n(a,q)&=\frac{1-aq^{2n}}{1-a}\sum_{k=0}^{n}\tau(n-k)
\frac{\poq{a}{n+k}}
{\poq{q}{n-k}}\beta_k(a,q).\label{invnew}
 \end{align}
\end{subequations}
\end{dy}
Closely related to  Bailey pairs is  the well-known  Bailey lemma originally due to  W. N. Bailey \cite[(3.1)]{bailey}. The reader may also consult \cite[Eq. (3.4.9)]{sla} or \cite[Eq. (3.33)] {andrews2} for further details.
\begin{yl}[The Bailey  lemma:  {\rm  \cite[(3.1)]{bailey}}]\label{baileylemmaold} For any Bailey pair  $(\alpha_n(a,q),\beta_n(a,q))$, there holds
\begin{align}
\sum_{n=0}^{\infty}\frac{\poq{x,y}{n}}
{\poq{aq/x,aq/y}{n}}\left(\frac{aq}{xy}\right)^n\alpha_n(a,q)\label{slater}\\
= \frac{\poq{aq,aq/xy}{\infty}} {\poq{aq/x,aq/y}{\infty}}
\sum_{n=0}^{\infty}\poq{x,y}{n}\left(\frac{aq}{xy}\right)^n\beta_n(a,q).\nonumber
 \end{align}
  \end{yl}

Particularly, for $a=1$ and $x,y\to\infty$, we have
 \begin{yl}[cf. {\rm \cite[Corollary 3]{ramanujan}}]\label{lemma4}For any Bailey pair  $(\alpha_n(1,q),\beta_n(1,q))$, there hold\begin{align}
\sum_{n=0}^{\infty} q^{n^2} \beta_n(1, q) & =\frac{1}{\poq{q}{\infty}} \sum_{n=0}^{\infty} q^{n^2} \alpha_n(1, q),\label{AAA1} \\
\sum_{n=0}^{\infty} q^{n^2}\left(-q ; q^2\right)_n \beta_n\left(1, q^2\right) & =\frac{1}{\psi(-q)} \sum_{n=0}^{\infty} q^{n^2} \alpha_n\left(1, q^2\right), \label{AAA2} \\
\sum_{n=0}^{\infty} q^{n(n+1) / 2}(-1;q)_n \beta_n(1, q) & =\frac{2}{\varphi(-q)} \sum_{n=0}^{\infty} \frac{q^{n(n+1) / 2}}{1+q^n} \alpha_n(1, q)\label{AAA3},\end{align}
where
$$ \psi(-q):=\frac{\left(q^2 ; q^2\right)_{\infty}}{\left(-q ; q^2\right)_{\infty}} ,~~ \varphi(-q):=\frac{(q ; q)_{\infty}}{(-q ; q)_{\infty}}.$$
\end{yl}

Furthermore,  we need a truncated Bailey transform for the Bailey pair $(\alpha_n(a,q),\beta_n(a,q))$ which is the case of \eqref{invnew} specialized with the Bailey pair  \cite[Ex. (13.105)]{ramanujan-book} due to J. Mc. Laughlin.
\begin{yl}[cf. \mbox{\cite[(13.105)]{ramanujan-book}}]For  any integer $m\geq 0$, it holds
\begin{align}
\tau(m)q^m\frac{(a q ; q)_m }{(q ; q)_m} \sum_{n=0}^m\left(q^{-m}, a q^{1+m} ; q\right)_n \beta_n(a, q)=\sum_{n=0}^m \alpha_n(a, q). \label{wangyee}
\end{align}
\end{yl}
For a good survey on Bailey pairs and various Bailey lemmas, as well as allied  extensions, the reader is referred to \cite{warnaar0} due to S. O. Warnaar. Especially noteworthy is that, as demonstrated by  G. E. Andrews and B. C. Berndt in \cite[\S 11.5]{ber}
and \cite[Chap. 5]{ber-1},   this  Bailey  lemma plays an important role in the study of Ramanujan's mock theta function identities.

In this  paper, we will show the following new Bailey pairs via a series of preliminaries. 
\begin{dl}\label{mainthm}
Two sequences $(\alpha_n(1,q),\beta_n(1,q))$  and $(\alpha_n^{(1)}(1,q),\beta_n^{(1)}(1,q))$ are two Bailey pairs with respect to $a=1$, where
\begin{align}\label{newbailey}
\left\{\begin{array}{l}
\alpha_n(1,q):=\tau(n)(1-q^n)G(n),\\
\\
\displaystyle
\beta_n(1,q):=\poq{q}{\infty}\frac{(q^n-1)\tau^2(n)}{\poq{q}{2n}}
\end{array}\right.
  \end{align}
and
\begin{align}\label{newbaileyg}
\left\{\begin{array}{l}
\alpha_n^{(1)}(1,q):=\tau(n)(1-q^{2n})H(n),\\
\\
\displaystyle
\beta_n^{(1)}(1,q):=\poq{q}{\infty}\frac{\tau^2(n){q^{-n}}}{\poq{q}{2n-1}}.
\end{array}\right.
\end{align}
Hereafter, we define
\begin{align}
G(n)&:=\sum_{k=-\infty}^\infty\,\tau^3(k)(1-q^k)q^{(1+n)k},\label{vvvvvv-0}\\
H(n)&:=\sum_{k=-\infty}^\infty\,\tau^3(k)q^{(3+n)k}.\label{vvvvvv-1}
\end{align}
\end{dl}
As applications of Theorem \ref{mainthm},  by combining these two Bailey pairs with the Bailey lemma and Lemma \ref{lemma4}, we will establish

\begin{dl}\label{dlsecond} For $xy\neq q$, there holds
\begin{align}
\sum_{n=0}^{\infty}\frac{\poq{x,y}{n}}{\poq{q}{2n}}
\left(\frac{1}{xy}\right)^n(q^n-1)q^{n^2}=(S_1-S_2)\frac {\poq{q/x,q/y}{\infty}}{\poq{q,q/xy}{\infty}}, \label{yyyyyy}
 \end{align}
where $S_1$  and  $S_2$ are, respectively, given by 
\begin{subequations}\label{xxxxxx}
 \begin{align}
S_1&:=\sum_{n=-\infty}^\infty\,\frac{\poq{x,y}{3n}}
{\poq{q/x,q/y}{3n}}\left(\frac{1}{xy}\right)^{3n}q^{3n^2+n}(q^n-1),\\
S_2&:=\sum_{n=-\infty}^\infty\,\frac{\poq{x,y}{3n+1}}
{\poq{q/x,q/y}{3n+1}}\left(\frac{1}{xy}\right)^{3n+1}q^{3n^2+4n+1}(1-q^{3n+1}).
 \end{align} 
 \end{subequations}
\end{dl}

\begin{dl}\label{dlthird}For $xy\neq q$, there holds
\begin{align}
\sum_{n=1}^{\infty}\frac{\poq{x,y}{n}}{\poq{q}{2n-1}}
\left(\frac{1}{xy}\right)^nq^{n^2-n}=(T_1-T_2)\frac {\poq{q/x,q/y}{\infty}}{\poq{q,q/xy}{\infty}}, \label{yyyyyyg}
 \end{align}
where $T_1$  and  $T_2$ are, respectively, given by 
\begin{subequations}\label{xxxxxxg}
 \begin{align}
 T_1&:=\sum_{n=-\infty}^\infty\,\frac{\poq{x,y}{3n+1}}
{\poq{q/x,q/y}{3n+1}}\left(\frac{1}{xy}\right)^{3n+1}q^{3 n^2+2n}\\
T_2&:=\sum_{n=-\infty}^\infty\,\frac{\poq{x,y}{3n+2}}
{\poq{q/x,q/y}{3n+2}}\left(\frac{1}{xy}\right)^{3n+2}q^{3 n^2+4n+1}.
 \end{align} 
 \end{subequations}
\end{dl}
As applications, we may deduce from these two theorems some $q$-identities modulo 15, 24, and 30. In particular, we may establish two $q$-identities very similar to  Ramanujan's lost notebook \cite[p. 4]{ramanujan-note}. The reader may consult  \cite[Entry. 12.3.4]{ber} for more details.
\begin{tl}\label{dlthird-one}For $x\neq 1$, the following identities hold:
\begin{align}
\sum_{n=-\infty}^\infty\,\frac{q^{3n^2+n}(q^n-1)}
{(1-xq^{3n})(1-q^{3n}/x)}&-\sum_{n=-\infty}^\infty\,\frac{q^{3n^2+4n+1}(1-q^{3n+1})}
{(1-xq^{3n+1})(1-q^{3n+1}/x)}\nonumber\\&=\frac{\poq{q}{\infty}^2} {\poq{x,1/x}{\infty}}
\sum_{n=0}^{\infty}\frac{\poq{x,1/x}{n}}{\poq{q}{2n}}
(q^n-1)q^{n^2}; \label{yyyyyy-one}\\
 \sum_{n=-\infty}^\infty\,\frac{q^{3n^2+2n}}
{(1-xq^{3n+1})(1-q^{3n+1}/x)}&-\sum_{n=-\infty}^\infty\,\frac{q^{3n^2+4n+1}}
{(1-xq^{3n+2})(1-q^{3n+2}/x)}\nonumber\\&=\frac{\poq{q}{\infty}^2} {\poq{xq,q/x}{\infty}}
\sum_{n=0}^{\infty}\frac{\poq{xq,q/x}{n}}{\poq{q}{2n+1}}
q^{n^2+n}. \label{yyyyyy-one}
 \end{align}
\end{tl}

\section{The proofs of the main results}
\setcounter{equation}{0}
For our purpose, we need to establish some preliminaries.  First of all, we need to show the following two generalizations of
  Euler's pentagonal number theorem. They  were first stated in \cite{wangjin} without any direct proof. As for Euler's pentagonal number theorem, we refer the reader to \cite{andrewseuler} for details. 

The following proof is due to S.O.Warnaar.
\begin{yl}[cf. \mbox{\cite[Eq. (4.16)/Eq. (4.17)]{wangjin}}]For any integer $n\geq 0$, it holds
\begin{align}
\sum_{k=-\infty}^\infty\,\tau^3(k)(q^{-k},q^k;q)_nq^k=\poq{q}{\infty}\tau^2(n);\label{cccccc}\\
\sum_{k=-\infty}^\infty\,\tau^3(k)(q^{-k},q^{k+1};q)_nq^{2k}
=\poq{q}{\infty}\tau^2(n){q^{n}}.\label{cccccc-new}
\end{align}
\end{yl}
\pf To show this lemma,   we begin with   (II.20) of \cite{10}, namely,
\begin{align*}\sum_{k=0}^\infty\frac{1-aq^{2k}}{1-a}
\frac{\poq{a,b,c,d}{k}}
 {\poq{q,aq/b,aq/c,aq/d}{k}}\bigg(\frac{aq}{bcd}\bigg)^k=\frac{(a q, a q / b c, a q / b d, a q / c d ; q)_{\infty}}{(a q / b, a q / c, a q / d, a q / b c d ; q)_{\infty}}
\end{align*}
and let $b,c,d\to \infty$. The result is
\begin{align*}
\sum_{k=0}^{\infty}(1-aq^{2k})\frac{\poq{a}{k}}
 {\poq{q}{k}}(aq)^k\tau^3(k)=\poq{a}{\infty}.
 \end{align*}
Next,  shift the summation index $k\to k-n$. Then we get
\begin{eqnarray}
\sum_{k=n}^{\infty}(1-aq^{2k-2n})\frac{\poq{a}{k-n}}
 {\poq{q}{k-n}}(aq)^{k-n}\tau^3(k-n)=\poq{a}{\infty}\label{2-baileyoperator-11}.
 \end{eqnarray}
Using the fact that $\tau(k-n)=\tau(k)q^{n(n-k)}/\tau(n)$ to simplify  \eqref{2-baileyoperator-11}
and replacing  $a$ with $aq^{2n}$,  
we have
\begin{eqnarray}
\sum_{k=n}^{\infty}(1-aq^{2k})\frac{\poq{aq^{2n}}{k-n}}
 {\poq{q}{k-n}}(aq^{1-n})^{k}\tau^3(k)=(aq^{1-n})^n\tau^3(n)\poq{aq^{2n}}{\infty}\label{2-baileyoperator-111}.
 \end{eqnarray}
Observe that
\begin{align*}
\poq{aq^{2n}}{k-n}=\frac{\poq{a}{k}\poq{aq^{k}}{n}}{\poq{a}{2n}},~~
\frac{1}{\poq{q}{k-n}}=\frac{\left(q^{-k}; q\right)_n}{(q ; q)_k}\frac{q^{nk}}{\tau(n)}.
\end{align*}
Substituting these expressions into \eqref{2-baileyoperator-111}, we find
\begin{eqnarray}
\sum_{k=0}^{\infty}(1-aq^{2k})\frac{\poq{a}{k}}{\poq{q}{k}}\left(q^{-k},aq^k; q\right)_n\tau^3(k)(aq)^{k}=(aq^{1-n})^n\tau^4(n)\poq{a}{\infty}
\label{2-baileyoperator-5}.
 \end{eqnarray}
In the following, we consider the cases $a=1$ and $a=q$ of \eqref{2-baileyoperator-5} .  The details are as follows.
\begin{description}
\item[(i) For $a=1$.] In this case, \eqref{2-baileyoperator-5} becomes to
\begin{align}
\sum_{k=0}^{\infty}(1+q^{k})\left(q^{-k},q^{k}; q\right)_n\tau^3(k)q^{k}=\tau^2(n)\poq{q}{\infty}.\label{rrr-1}
\end{align}
Note that
\begin{align*}
\sum_{k=0}^{\infty}\left(q^{-k},q^{k}; q\right)_n\tau^3(k)q^{2k}=\sum_{k=-\infty}^{0}\poq{q^{-k},q^{k}}{n}\tau^3(k)q^k.
\end{align*}
Then \eqref{rrr-1} turns out to be \eqref{cccccc}.
\item[(ii) For $a=q$.] As such, \eqref{2-baileyoperator-5} reduces to
\begin{align}
&\sum_{k=0}^{\infty}(1-q^{2k+1})\left(q^{-k},q^{k+1}; q\right)_n\tau^3(k)q^{2k}=q^n\tau^2(n)\poq{q}{\infty}.\label{rrr-2}
\end{align}
It is clear that
\begin{align*}
-\sum_{k=0}^{\infty}\left(q^{-k},q^{k+1}; q\right)_n\tau^3(k)q^{4k+1}\stackrel{k\to -k-1}{===}\sum_{k=-\infty}^{-1}\left(q^{-k},q^{k+1}; q\right)_n\tau^3(k)q^{2k},
\end{align*}
which leads us to \eqref{cccccc-new} at once.
\end{description}
Summing up, the lemma is proved.
   \qed

\begin{yl}\label{prop-h}Let $G(n)$ and $H(n)$ be defined respectively by \eqref{vvvvvv-0} and \eqref{vvvvvv-1}. Then we have
\begin{align}
\frac{G(n)}{\poq{q}{\infty}}&=\left\{\begin{array}{ll}
(-1)^L(q^L-1)q^{-3L^2/2-L/2},&\mbox{if}~~n=3L;\\
(-1)^Lq^{-3L^2/2-L/2},&\mbox{if}~~n=3L+1;\\
(-1)^Lq^{-3L^2/2-5L/2-1},&\mbox{if}~~n=3L+2.
\end{array}\right.\label{vvvvvv}\\
\mbox{and}\qquad\qquad&\nonumber\\
\frac{H(n)}{\poq{q}{\infty}}&=\left\{\begin{array}{ll}
0,&\mbox{if}~~n=3L;\\
(-1)^{L+1}q^{-3L^2/2-5L/2-1},&\mbox{if}~~n=3L+1;\\
(-1)^{L+1}q^{-3L^2/2-7L/2-2},&\mbox{if}~~n=3L+2.
\end{array}\right.\label{vvvvvvg}
\end{align}
\end{yl}
\pf  To show \eqref{vvvvvv} and \eqref{vvvvvvg}, it only needs to apply Jacobi's triple product identity \eqref{triple} to the sums
\begin{align*}
G(n)&=\sum_{k=-\infty}^\infty\,\tau^3(k)q^{(1+n)k}-\sum_{k=-\infty}^\infty\,\tau^3(k)q^{(2+n)k}\\
&=(q^{1+n},q^{2-n},q^3;q^3)_\infty-(q^{2+n},q^{1-n},q^3;q^3)_\infty,\\
H(n)&=\sum_{k=-\infty}^\infty\,\tau^3(k)q^{(3+n)k}=(q^{3+n},q^{-n},q^3;q^3)_\infty.
\end{align*}
Next, according as   $n$ modulo $3$,  we  evaluate these two sums as desired.
\qed

Making use of  Lemma \ref{prop-h}, it is easy to verify  that 
\begin{prop}\label{GHprops}
$G(0)=H(0)=0$ and
\begin{align}G(-n)=-G(n);~~ H(-n)=-q^nH(n).
\end{align}
\end{prop}
Besides, both $G(n)$ and $H(n)$ satisfy the following $q$-identities which are the key ingredients of  our discussions.
\begin{yl}Let $G(n)$ and $H(n)$ be given as above. Then, for any integer $n\geq 0$, it holds 
\begin{align}
\sum_{k=0}^n\genfrac{[}{]}{0pt}{}{2n}{k+n}_{q}\tau(k)(1-q^k)G(k)&=\poq{q}{\infty}(q^n-1)\tau^2(n),\label{dddddd}\\
\sum_{k=0}^n\genfrac{[}{]}{0pt}{}{2n}{k+n}_{q}\frac{1-q^{2k}}{1-q^{2n}}\tau(k)H(k)
&=\poq{q}{\infty}\tau^2(n){q^{-n}}.\label{ddddddg}
\end{align}
\end{yl}
\pf (i) In order to establish \eqref{dddddd},  we need to start from  \eqref{cccccc}  and   expand
\begin{align*}
(x,1/x;q)_{n}=\prod_{i=0}^{n-1}(1-xq^i)\prod_{i=0}^{n-1}(1-q^i/x)=\tau(n)(1-x)
\poq{xq^{1-n}}{2n-1}x^{-n},
\end{align*}from which it follows
\begin{align*}
(q^{-k},q^k;q)_{n}=\tau(n)q^{-nk}(1-q^k)\poq{q^{1-n+k}}{2n-1}.
\end{align*}
Substituting this expression into \eqref{cccccc} and then expanding $\poq{q^{1-n+k}}{2n-1}$ by the q-binomial theorem \cite[(II.3)]{10}, we get
\begin{align*}
\poq{q}{\infty}\tau(n)&=\sum_{k=-\infty}^\infty\,\tau^3(k)q^{k-nk}(1-q^k)\sum_{i=0}^{2n-1}\genfrac{[}{]}{0pt}{}{2n-1}{i}_{q}\tau(i)(q^{1-n+k})^i\\
&=\sum_{i=0}^{2n-1}\genfrac{[}{]}{0pt}{}{2n-1}{i}_{q}\tau(i)q^i\sum_{k=-\infty}^\infty\,\tau_3(k)(1-q^k)q^{k-nk+ki-ni}.
\end{align*}
Now, shift the summation index $i$ to $i+n$. Then we get
\begin{align*}
\poq{q}{\infty}\tau(n)=\tau(n)q^{-n^2+n}\sum_{i=-n}^{n-1}\genfrac{[}{]}{0pt}{}{2n-1}{i+n}_{q}\tau(i)q^{i}\sum_{k=-\infty}^\infty\,\tau^3(k)(1-q^k)q^{k+ki},
\end{align*}
which, after simplified by  Proposition \ref{GHprops}, reduces to
\begin{align*}
&\poq{q}{\infty}\tau^2(n)
=\sum_{i=-n}^{n-1}\genfrac{[}{]}{0pt}{}{2n-1}{i+n}_{q}\tau(i)q^{i}G(i)\\
&\quad=\bigg\{\sum_{i=1}^{n}+\sum_{i=-n}^{-1}\bigg\}\genfrac{[}{]}{0pt}{}{2n-1}{i+n}_{q}\tau(i)q^{i}G(i)=\sum_{i=1}^{n}\bigg\{\genfrac{[}{]}{0pt}{}{2n-1}{i+n}_{q}q^{i}-\genfrac{[}{]}{0pt}{}{2n-1}{n-i}_{q}\bigg\}\tau(i)G(i).
\end{align*}
Observe  that
\[\genfrac{[}{]}{0pt}{}{2n-1}{i+n}_{q}q^{i}-\genfrac{[}{]}{0pt}{}{2n-1}{n-i}_{q}=-\genfrac{[}{]}{0pt}{}{2n}{i+n}_{q}\frac{1-q^i}{1-q^n},\]
yielding the desired identity \eqref{dddddd}. 

(ii) To establish \eqref{ddddddg},  as have done as above, we first  expand
\begin{align*}
(q^{-k},q^{k+1};q)_{n}&=\prod_{i=0}^{n-1}(1-q^{-k+i})\prod_{i=0}^{n-1}(1-q^{k+1+i})\\
&=\tau(n)q^{-nk}\prod_{i=k-n+1}^k(1-q^i)\prod_{i=k+1}^{n+k}(1-q^{i})
=\tau(n)q^{-nk}(q^{1-n+k};q)_{2n}.
\end{align*}
Substituting this expression into \eqref{cccccc-new} and then expanding $(q^{1-n+k};q)_{2n}$ by the q-binomial theorem \cite[(II.3)]{10}, we get
\begin{align*}
\poq{q}{\infty}\tau(n){q^n}&=\sum_{k=-\infty}^\infty\,\tau^3(k)q^{2k-nk}\sum_{i=0}^{2n}\genfrac{[}{]}{0pt}{}{2n}{i}_{q}\tau(i)(q^{1-n+k})^i\\
&=\sum_{i=0}^{2n}\genfrac{[}{]}{0pt}{}{2n}{i}_{q}\tau(i)q^i\sum_{k=-\infty}^\infty\,\tau^3(k)q^{2k-nk+ki-ni}.
\end{align*}
Now, shift the summation index $i$ to $i+n+1$ and simplify further. Then we get
\begin{align*}
\poq{q}{\infty}\tau(n){q^n}=\tau(n)q^{-n^2+n}\sum_{i=-n-1}^{n-1}
\genfrac{[}{]}{0pt}{}{2n}{i+n+1}_{q}\tau(i+1)q^{i+1}\sum_{k=-\infty}^\infty\,\tau^3(k)q^{3k+ki}.
\end{align*}
As above,  we simplify this identity by Proposition \ref{GHprops}, obtaining 
\begin{align*}
\poq{q}{\infty}\tau^2(n){q^n}&=\bigg\{\sum_{i=0}^{n+1}+\sum_{i=-n-1}^{0}\bigg\}
\genfrac{[}{]}{0pt}{}{2n}{i+n+1}_{q}\tau(i+1)q^{i+1}H(i)\\
&=q\sum_{i=0}^{n+1}\bigg\{\genfrac{[}{]}{0pt}{}{2n}{-i+n+1}_{q}
-\genfrac{[}{]}{0pt}{}{2n}{i+n+1}_{q}q^{2i}\bigg\}\tau(i)H(i).
\end{align*}
Now shift the summation index $n$ to $n-1$. We get
\begin{align*}
\poq{q}{\infty}\tau^2(n-1){q^{n-1}}=q\sum_{i=0}^{n}\bigg\{\genfrac{[}{]}{0pt}{}{2n-2}{-i+n}_{q}
-\genfrac{[}{]}{0pt}{}{2n-2}{i+n}_{q}q^{2i}\bigg\}\tau(i)H(i).
\end{align*}
It is easy to check
\[\genfrac{[}{]}{0pt}{}{2n-2}{-i+n}_{q}
-\genfrac{[}{]}{0pt}{}{2n-2}{i+n}_{q}q^{2i}=\genfrac{[}{]}{0pt}{}{2n}{i+n}_{q}\frac{1-q^{2i}}{1-q^{2n}},\]
yielding the desired identity \eqref{ddddddg}.  The Lemma is proved.
\qed

Now, we are in a good position to show Theorem \ref{mainthm}. 

\begin{pot444}  Obviously, \eqref{dddddd} can be reformulated as the form
\begin{align}
\sum_{k=0}^n\frac{\tau(k)(1-q^k)G(k)}{\poq{q}{n-k}\poq{q}{n+k}}=\poq{q}{\infty}\frac{(q^n-1)\tau^2(n)}{\poq{q}{2n}}.\label{dddddd-bailey}
\end{align}
while \eqref{ddddddg} may be reformulated as
\begin{align}
\sum_{k=0}^n\frac{\tau(k)(1-q^{2k})H(k)}{\poq{q}{n-k}\poq{q}{n+k}}
=\poq{q}{\infty}\frac{\tau^2(n){q^{-n}}}{\poq{q}{2n-1}}.\label{ddddddg-bailey}
\end{align}
Obviously,  in light of Definition \ref{defnew-ooo}, we have the conclusion.\qed
\end{pot444}
\begin{remark} As S. O. Warnaar pointed out,   the Bailey pair  $(\alpha_n(1,q),\beta_n(1,q))$  turns out to be the $A(5)-A(7)$, where $A(5)$ and $A(7)$ are Bailey pairs  in the paper \cite{Slater} due to L. J. Slater.
\end{remark}Furthermore, on considering the  inverse forms of \eqref{dddddd-bailey} and \eqref{ddddddg-bailey}, we have 
\begin{prop} Let $G(n)$ and $H(n)$ be given by \eqref{vvvvvv-0} and \eqref{vvvvvv-1}, respectively. Then
\begin{align}
G(n)&=\poq{q}{\infty}\sum_{k=1}^n(-1)^{k+1} \frac{\poq{q^{2k}}{n-k}}{\poq{q}{n-k}}\frac{1+q^n}{1+q^k}q^{3k^2/2-k/2-nk},\label{qid18}\\
H(n)&=\poq{q}{\infty}\sum_{k=1}^n(-1)^{k} \frac{\poq{q^{2k}}{n-k}}{\poq{q}{n-k}}q^{3k^2/2-{3}k/2-nk}.\label{qid18-18}
\end{align}
\end{prop}

Having established such two pairs of Bailey pairs, we are able  to show  the other main results, namely, Theorems \ref{dlsecond}-\ref{dlthird}, as well as Corollary \ref{dlthird-one}.
\begin{pot555}
 It suffices to  set  $a=1$ and specificize the Bailey lemma \eqref{slater} with the  Bailey pair $(\alpha_n(1,q),\beta_n(1,q))$. Then \eqref{slater} reduces to
\begin{align}
\sum_{n=0}^{\infty}\frac{\poq{x,y}{n}}
{\poq{q/x,q/y}{n}}\left(\frac{q}{xy}\right)^n\tau(n)(1-q^n)
G(n)\label{slater-111}\\
= \frac{\poq{q,q/xy}{\infty}} {\poq{q/x,q/y}{\infty}}
\sum_{n=0}^{\infty}\poq{x,y}{n}\left(\frac{q}{xy}\right)^n\poq{q}{\infty}\frac{(q^n-1)\tau^2(n)}{\poq{q}{2n}}.\nonumber
 \end{align}
 Next, in view of the fact that $G(-n)=-G(n)$, we may reformulate it as a  bilateral series. It is because that
 \begin{align*}
\mbox{LHS of \eqref{slater-111}}&=\sum_{n=0}^{\infty}\frac{\poq{x,y}{n}}
{\poq{q/x,q/y}{n}}\left(\frac{q}{xy}\right)^n\tau(n)G(n)\\
&+\sum_{n=0}^{\infty}\frac{\poq{x,y}{n}}
{\poq{q/x,q/y}{n}}\left(\frac{q}{xy}\right)^n
\tau(n)q^nG(-n)
\end{align*}
 while the second series
\begin{align*}
\sum_{n=0}^{\infty}\frac{\poq{x,y}{n}}
{\poq{q/x,q/y}{n}}\left(\frac{q}{xy}\right)^n\tau(n)q^nG(-n)
&=\sum_{n=-\infty}^{0}\frac{\poq{x,y}{-n}}
{\poq{q/x,q/y}{-n}}\left(\frac{q}{xy}\right)^{-n}\tau(-n)q^{-n}G(n)\\
&=\sum_{n=-\infty}^{0}\frac{\poq{x,y}{n}}
{\poq{q/x,q/y}{n}}\left(\frac{q}{xy}\right)^n
\tau(n)G(n),
\end{align*}
thus the left-hand side of \eqref{slater-111} indeed forms a bilateral series 
\begin{align*} 
\sum_{n=-\infty}^\infty\,\frac{\poq{x,y}{n}}
{\poq{q/x,q/y}{n}}\left(-\frac{1}{xy}\right)^nq^{n(n+1)/2}G(n).
\end{align*}
 Therefore
\begin{align}
& \frac{(q, q / x y ; q)_{\infty}}{(q / x, q / y ; q)_{\infty}} \sum_{n=0}^{\infty}\left(q^n-1\right) q^{n^2} \frac{(x, y ; q)_n}{(q ; q)_{2 n}}\left(\frac{1}{x y}\right)^n \nonumber\\
& =\sum_{n=-\infty}^\infty\,\frac{(x, y ; q)_n}{(q / x, q / y ; q)_n}\left(-\frac{1}{x y}\right)^n q^{n(n+1) / 2} \frac{G(n)}{(q ; q)_{\infty}}.\label{slater-111p-222}
\end{align}
Next, we insert the expression of $G(n)/\poq{q}{\infty}$ given by \eqref{vvvvvv} into \eqref{slater-111p-222} and simplify further. The detail is as follows:
\begin{align*}
\frac{\poq{q,q/xy}{\infty}} {\poq{q/x,q/y}{\infty}}
&\sum_{n=0}^{\infty}(q^n-1)q^{n^2}\frac{\poq{x,y}{n}}{\poq{q}{2n}}
\left(\frac{1}{xy}\right)^n\\
&=\sum_{n=-\infty}^\infty\,\frac{\poq{x,y}{3n}}
{\poq{q/x,q/y}{3n}}\left(-\frac{1}{xy}\right)^{3n}q^{3n(3n+1)/2}\frac{G(3n)}{\poq{q}{\infty}}\\
&+\sum_{n=-\infty}^\infty\,\frac{\poq{x,y}{3n+1}}
{\poq{q/x,q/y}{3n+1}}\left(-\frac{1}{xy}\right)^{3n+1}q^{(3n+1)(3n+2)/2}\frac{G(3n+1)}{\poq{q}{\infty}}\\
&+\sum_{n=-\infty}^\infty\,\frac{\poq{x,y}{3n+2}}
{\poq{q/x,q/y}{3n+2}}\left(-\frac{1}{xy}\right)^{3n+2}q^{(3n+2)(3n+3)/2}\frac{G(3n+2)}{\poq{q}{\infty}}
 \\
 &=S_1-S_2, \end{align*}
where $S_i~(1\leq i\leq 2)$ are given by \eqref{xxxxxx}. The theorem is confirmed.
  \qed
\end{pot555}\
As performing above, we specialize the Bailey lemma to Theorem \ref{dlthird}. The following is the full proof.
\begin{pot666}  Given the Bailey pair $(\alpha_n^{(1)}(1,q),\beta_n^{(1)}(1,q))$,    the Bailey lemma \eqref{slater} with $a=1$ reduces to
\begin{align}
\sum_{n=0}^{\infty}\frac{\poq{x,y}{n}}
{\poq{q/x,q/y}{n}}\left(\frac{q}{xy}\right)^n\tau(n)(1-q^{2n})
H(n)\label{slater-222}\\
= \frac{\poq{q,q/xy}{\infty}} {\poq{q/x,q/y}{\infty}}
\sum_{n=1}^{\infty}\poq{x,y}{n}\left(\frac{q}{xy}\right)^n\poq{q}{\infty}\frac{\tau^2(n)q^{-n}}{\poq{q}{2n-1}}.\nonumber
 \end{align}
 Next, by  using the fact that $H(-n)=-q^nH(n)$, we can still  reformulate it as a  bilateral series. Since
 \begin{align*}
\mbox{LHS of \eqref{slater-222}}&=\sum_{n=0}^{\infty}\frac{\poq{x,y}{n}}
{\poq{q/x,q/y}{n}}\left(\frac{q}{xy}\right)^n\tau(n)H(n)\\
&+\sum_{n=0}^{\infty}\frac{\poq{x,y}{n}}
{\poq{q/x,q/y}{n}}\left(\frac{q}{xy}\right)^n
\tau(n)q^nH(-n)
\end{align*}
 while the second series
\begin{align*}
\sum_{n=0}^{\infty}\frac{\poq{x,y}{n}}
{\poq{q/x,q/y}{n}}\left(\frac{q}{xy}\right)^n\tau(n)q^nH(-n)
&=\sum_{n=-\infty}^{0}\frac{\poq{x,y}{-n}}
{\poq{q/x,q/y}{-n}}\left(\frac{q}{xy}\right)^{-n}\tau(-n)q^{-n}H(n)\\
&=\sum_{n=-\infty}^{0}\frac{\poq{x,y}{n}}
{\poq{q/x,q/y}{n}}\left(\frac{q}{xy}\right)^n
\tau(n)H(n),
\end{align*}
thus the left-hand side of \eqref{slater-222} is still  a bilateral series 
\begin{align*} 
\sum_{n=-\infty}^\infty\,\frac{\poq{x,y}{n}}
{\poq{q/x,q/y}{n}}\left(-\frac{1}{xy}\right)^nq^{n(n+1)/2}H(n).
 \end{align*}
 Therefore
\begin{align}
& \frac{(q, q / x y ; q)_{\infty}}{(q / x, q / y ; q)_{\infty}} \sum_{n=1}^{\infty} q^{n^2-n} \frac{(x, y ; q)_n}{(q ; q)_{2 n-1}}\left(\frac{1}{x y}\right)^n =\sum_{n=-\infty}^\infty\,\frac{(x, y ; q)_n}{(q / x, q / y ; q)_n}\left(-\frac{1}{x y}\right)^n q^{n(n+1) / 2} \frac{H(n)}{(q ; q)_{\infty}}.\label{slater-111p-222g}
\end{align}
Next, we insert the expression of $H(n)/\poq{q}{\infty}$ given by \eqref{vvvvvvg} into \eqref{slater-111p-222g} and simplify further. The detail is as follows:
\begin{align*}
&\quad\frac{\poq{q,q/xy}{\infty}} {\poq{q/x,q/y}{\infty}}
\sum_{n=1}^{\infty}q^{n^2-n}\frac{\poq{x,y}{n}}{\poq{q}{2n-1}}
\left(\frac{1}{xy}\right)^n\\
&=\sum_{n=-\infty}^\infty\,\frac{\poq{x,y}{3n+1}}
{\poq{q/x,q/y}{3n+1}}\left(-\frac{1}{xy}\right)^{3n+1}q^{(3n+1)(3n+2)/2}\frac{H(3n+1)}{\poq{q}{\infty}}\\
&+\sum_{n=-\infty}^\infty\,\frac{\poq{x,y}{3n+2}}
{\poq{q/x,q/y}{3n+2}}\left(-\frac{1}{xy}\right)^{3n+2}q^{(3n+2)(3n+3)/2}\frac{H(3n+2)}{\poq{q}{\infty}}
 \\
 &=T_1-T_2, \end{align*}
where $T_i ~ (i=1,2)$ are given by \eqref{xxxxxxg}. The theorem is confirmed.\qed
\end{pot666}
As specific cases of  Theorems \ref{dlsecond} and \ref{dlthird}, we can show Corollary \ref{dlthird-one} in one-line argument.
\begin{pot777} It suffices to set $y=1/x$ in both Theorem \ref{dlsecond} and Theorem \ref{dlthird}. The rest are routine manipulations.
\end{pot777}

\section{Applications: $q$-identities based on Bailey pairs}
\setcounter{equation}{0}
\subsection{$q$-Series identities based on  $(\alpha_n(1,q),\beta_n(1,q))$}
Having established such a Bailey pair $(\alpha_n(1,q),\beta_n(1,q))$ given by \eqref{newbailey}, we proceed to establish concrete $q$-identities. Indeed, as good justification, we may derive from \eqref{yyyyyy}  that 
 \begin{tl} [Mod 30 identity]
\begin{align}
\sum_{n=0}^{\infty}\frac{\poqq{-q}{n}}{\poqq{q^2}{2n-1}}q^{3n^2-2n}
&=\frac {\poqq{-q}{\infty}}{\poqq{q^2}{\infty}}\bigg((-q^{11},-q^{19},q^{30};q^{30})_\infty-(-q^{13},-q^{17},q^{30};q^{30})_\infty\label{yyyyyy-new}\\
&+q(-q^{7},- q^{23},q^{30};q^{30})_\infty-q^3(-q,-q^{29},q^{30};q^{30})_\infty\bigg). \nonumber
\end{align} 
\end{tl}
\pf It suffices to replace $q$ with $q^2$  in \eqref{yyyyyy} first.  Then
 \begin{align}
\sum_{n=0}^{\infty}(q^{2n}-1)q^{2n^2}\frac{\poqq{x,y}{n}}{\poqq{q^2}{2n}}
\left(\frac{1}{xy}\right)^n=(S_1-S_2)\frac {\poqq{q^2/x,q^2/y}{\infty}}{\poqq{q^2,q^2/xy}{\infty}}. \label{yyyyyy-two}
 \end{align}
 And then set $y\to \infty$ alone. Then \eqref{yyyyyy-two} reduces to
  \begin{align}
\sum_{n=0}^{\infty}(q^{2n}-1)q^{3n^2-n}\frac{\poqq{x}{n}}{\poqq{q^2}{2n}}\left(-\frac{1}{x}\right)^n
=(S_1-S_2)\frac {\poqq{q^2/x}{\infty}}{\poqq{q^2}{\infty}}, \label{temp-1}
 \end{align}
 In this form, we further set $x=-q$ and evaluate each $S_i$  in \eqref{xxxxxx} via Jacobi's triple product identity  \eqref{triple}. We obtain
 \begin{align*}
&S_1=\sum_{n=-\infty}^\infty\,q^{15n^2-2n}-\sum_{n=-\infty}^\infty\,q^{15n^2-4n}=(-q^{13},- q^{17},q^{30};q^{30})_\infty-(-q^{11},- q^{19},q^{30};q^{30})_\infty,\\
&S_2=\sum_{n=-\infty}^\infty\,q^{15n^2+8n+1}-\sum_{n=-\infty}^\infty\,q^{15n^2+16n+4}=q(-q^{7},- q^{23},q^{30};q^{30})_\infty-q^3(-q,-q^{29},q^{30};q^{30})_\infty.
 \end{align*} 
A direct substitution of these results into \eqref{temp-1} leads us to \eqref{yyyyyy-new}.
  \qed
  
\begin{lz} Actually,  once applying the Bailey pair $(\alpha_n(1,q^2),\beta_n(1,q^2))$ to \eqref{AAA2} directly, we  can also deduce \eqref{yyyyyy-new}:
\begin{align}
 \sum_{n=0}^{\infty}\frac{(-q ; q^2)_n}{\poqq{q^2}{2n}} (q^{2n}-1) q^{3n^2-2n} &=\frac{(q^{30};q^{30})_\infty}{\psi(-q)}\bigg((-q^{17},-q^{13};q^{30})_\infty-(-q^{11},-q^{19};q^{30})_\infty\nonumber\\
&-q(-q^{7},-q^{23};q^{30})_\infty+q^3(-q,-q^{29};q^{30})_\infty\bigg).\label{188}
\end{align}
It should be pointed out that \eqref{188} is just the difference of Eq. (3.30.5) and Eq. (3.30.6) of 
\cite{ramanujan}.
\end{lz}

\begin{lz}
From \eqref{yyyyyy}   it follows  a $q$-identity of Rogers-Ramanujan type modulo $24$, as below: 
 \begin{align}
\poq{q}{\infty}\sum_{n=1}^{\infty} \frac{q^{2n^2-n}}{(1+q^n)\poq{q}{2n-1}}
&=(-q^{10},-q^{14},q^{24};q^{24})_\infty-(-q^{11},-q^{13},q^{24};q^{24})_\infty\nonumber\\
&+q(-q^5,-q^{19},q^{24};q^{24})_\infty-q^2(-q^{2},-q^{22},q^{24};q^{24})_\infty.\label{ooo-1}
\end{align}
This result  turns out to be the difference of Eq. (3.24.3) and Eq. (3.24.1) of \cite{ramanujan}. 
In fact, upon taking the limitation $x,y\to \infty$ of  \eqref{yyyyyy} simultaneously, we easily find 
\begin{align*}
\sum_{n=0}^{\infty}q^{2n^2-n}\frac{q^n-1}{\poq{q}{2n}}
=\frac {S_1-S_2}{\poq{q}{\infty}}. 
 \end{align*}
In such case, we find
 \begin{align*}
S_1&=\sum_{n=-\infty}^\infty\,q^{12n^2-n}-\sum_{n=-\infty}^\infty\,q^{12n^2-2n},\\
S_2&=\sum_{n=-\infty}^\infty\,q^{12n^2+7n+1}-\sum_{n=-\infty}^\infty\,q^{12n^2+14n+4}.
 \end{align*} 
By virtue of Jacobi's triple product identity \eqref{triple}, we are able to evaluate each $S_i$, thereby reducing \eqref{yyyyyy} to \eqref{ooo-1}.
\end{lz}

Upon substituting the Bailey pair $(\alpha_n(1,q),\beta_n(1,q))$ into \eqref{AAA3}, we may obtain
\begin{tl}[Mod 15 identity]\label{NNNNNN}
\begin{align}
\sum_{n=0}^{\infty}\frac{(-q;q)_{n}}{\poq{q}{2n+1}} \frac{q^{3n^2/2+5n/2+1}}{1+q^{n+1}}&+\frac{2}{\varphi(-q)}\sum_{n=-\infty}^{\infty}\bigg(\frac{q^{15n^2/2+n/ 2}}{1+q^{3n}}-\frac{q^{15n^2/2+11n/ 2+1}}{1+q^{3n+1}}\bigg)\nonumber
\\
&=\frac{1}{\varphi(-q)}\bigg((-q^7,-q^8,q^{15};q^{15})_{\infty}-q(-q^2,-q^{13},q^{15};q^{15})_{\infty}\bigg).\label{220}
\end{align}
\end{tl}
\pf It suffices to insert $(\alpha_n(1,q),\beta_n(1,q))$ into \eqref{AAA3}, yielding
\begin{align*}
\poq{q}{\infty}&\sum_{n=0}^{\infty}(q^n-1) q^{3n^2/2-n/2} \frac{(-1;q)_n}{\poq{q}{2n}}=\frac{2}{\varphi(-q)} \sum_{n=0}^{\infty} \frac{q^{n(n+1) / 2}}{1+q^n} \tau(n)(1-q^n)G(n)\\
&=\frac{2\poq{q}{\infty}}{\varphi(-q)}\sum_{n=0}^{\infty} \frac{q^{3n(3n+1) / 2}}{1+q^{3n}}\tau(3n)(1-q^{3n})(-1)^n(q^n-1)q^{-3n^2/2-n/2}\\
&+\frac{2\poq{q}{\infty}}{\varphi(-q)}\sum_{n=0}^{\infty}\frac{q^{(3n+1)(3n+2) / 2}}{1+q^{3n+1}}\tau(3n+1)(1-q^{3n+1})(-1)^nq^{-3n^2/2-n/2}\\
&+\frac{2\poq{q}{\infty}}{\varphi(-q)}\sum_{n=0}^{\infty} \frac{q^{3(n+2)(3n+3) / 2}}{1+q^{3n+2}}\tau(3n+2)(1-q^{3n+2})(-1)^nq^{-3n^2/2-5n/2-1}.
\end{align*}
After a bit simplification, we arrive at
\begin{align*}
&\sum_{n=0}^{\infty}(q^n-1) q^{3n^2/2-n/2} \frac{(-1;q)_n}{\poq{q}{2n}}\\
&=\frac{2}{\varphi(-q)}\left(\sum_{n=0}^{\infty} \frac{1-q^{3n}}{1+q^{3n}}q^{15n^2/2+n/ 2}-\sum_{n=0}^{\infty} \frac{1-q^{3n}}{1+q^{3n}}q^{15n^2/2-n/ 2}\right)\\
&-\frac{2}{\varphi(-q)}\sum_{n=0}^{\infty}
\frac{1-q^{3n+1}}{1+q^{3n+1}}q^{15n^2/2+11n/ 2+1}+\frac{2}{\varphi(-q)}\sum_{n=0}^{\infty} \frac{1-q^{3n+2}}{1+q^{3n+2}}q^{15n^2/2+19n/ 2+3}\\
&=\frac{2}{\varphi(-q)}\sum_{n=-\infty}^{\infty} \frac{1-q^{3n}}{1+q^{3n}}q^{15n^2/2+n/ 2}-\frac{2}{\varphi(-q)}\sum_{n=-\infty}^{\infty}
\frac{1-q^{3n+1}}{1+q^{3n+1}}q^{15n^2/2+11n/ 2+1}\\
&=\frac{4}{\varphi(-q)}\sum_{n=-\infty}^{\infty}\bigg(\frac{q^{15n^2/2+n/ 2}}{1+q^{3n}}-\frac{q^{15n^2/2+11n/ 2+1}}{1+q^{3n+1}}\bigg)-\frac{2}{\varphi(-q)}\sum_{n=-\infty}^{\infty}\bigg(q^{15n^2/2+n/ 2}-q^{15n^2/2+11n/ 2+1}\bigg).
\end{align*}
This is, after a bit simplification and one use of Jacobi's triple product identity,  in agreement with \eqref{220}.
\qed

We end this part with  the following summation formula for $q$-series, which is similar to but different from \eqref{qid18}.
 \begin{dl}For any integer $m\geq 0$, it holds
 \begin{align}\sum_{n=0}^m\tau^2(n)\tau(n-m) \frac{(q^n-1)\left(q^{m-n+1}; q\right)_{2n}}{\poq{q}{2n}}=\left\{\begin{array}{ll}
(1-q^L)q^{3L^2+L},&\mbox{if}~~m=3L;\\
(q^{2L+1}-1)q^{3L^2+2L},&\mbox{if}~~m=3L+1;\\
(1-q^{L+1})q^{3L^2+4L+1},&\mbox{if}~~m=3L+2.
\end{array}\right.\label{wangyee-1}
\end{align}
 \end{dl}
 \pf It suffices to set $a=1$ in \eqref{wangyee} and substitute $(\alpha_n(1, q),\beta_n(1, q))$ into. Then we have
 \begin{align}
\tau(m)q^{m}\sum_{n=0}^m\left(q^{-m}, q^{1+m} ; q\right)_n \frac{(q^n-1)\tau^2(n)}{\poq{q}{2n}}=F(m),\label{wangyee-2}
\end{align}
where, for convenience we  write 
\[F(m):=\sum_{k=0}^m \tau(k)(1-q^k)\frac{G(k)}{\poq{q}{\infty}}.\]
Using the basic relations that $\tau(n-m)=\tau(n)\tau(-m)q^{-mn}$ and 
\begin{align}
\left(q^{-m}, q^{1+m} ; q\right)_n=\poq{q^{m-n+1}}{2n}\tau(n)q^{-mn},\label{basictworelas}\end{align}
 we further simplify
\begin{align*}
\mbox{LHS of \eqref{wangyee-2}}=
\sum_{n=0}^m(q^n-1)\tau^2(n)\tau(n-m) \frac{\left(q^{m-n+1}; q\right)_{2n}}{\poq{q}{2n}}.
\end{align*}
Thanks to \eqref{vvvvvv}, we are able to  evaluate $F(m)$ on the right-hand side of  \eqref{wangyee-2} according as the cases $m\pmod 3$. First of all, let $m=3L$. We compute directly 
\begin{align*}
F(3L)&=
\sum_{k=0}^{L}\tau(3k)(1-q^{3k})\frac{G(3k)}{\poq{q}{\infty}}
+\sum_{k=0}^{L-1}\tau(3k+1)(1-q^{3k+1})\frac{G(3k+1)}{\poq{q}{\infty}}\\
&+\sum_{k=0}^{L-1}\tau(3k+2)(1-q^{3k+2})\frac{G(3k+2)}{\poq{q}{\infty}}\\
&=\sum_{k=0}^{L}(1-q^{3k})(q^k-1)q^{3k^2-2k}-\sum_{k=0}^{L-1}(1-q^{3k+1})q^{3k^2+k}\\
&+\sum_{k=0}^{L-1}(1-q^{3k+2})q^{3k^2+2k}=(1-q^L)q^{3L^2+L}.
\end{align*}
Next, when $m=3L+1$  and $m=3L+2$, it is clear that 
\begin{align*}
F(3L+1)&=F(3L)-(1-q^{3L+1})q^{3L^2+L}=(q^{2L+1}-1)q^{3L^2+2L},\\
F(3L+2)&=F(3L+1)+(1-q^{3L+2})q^{3L^2+2L}=(1-q^{L+1})q^{3L^2+4L+1}.
\end{align*}
Simplify \eqref{wangyee-2} by these expressions. We finally get \eqref{wangyee-1}.
 \qed
 
\subsection{$q$-Series identities based on the second Bailey pair $(\alpha_n^{(1)}(1,q),\beta_n^{(1)}(1,q))$}

In this section, we will focus on   concrete $q$-identities based on the Bailey pair $(\alpha_n^{(1)}(1,q),\beta_n^{(1)}(1,q))$ given by \eqref{newbaileyg}. As is expected,  we now derive  one  $q$-identity modulo 30 from \eqref{yyyyyyg}, which is the same as Eq. (3.30.7) of \cite{ramanujan}. 
 \begin{tl} [Mod 30 identity]
\begin{align}
&  \sum_{n=0}^{\infty}\frac{\poqq{-q}{n+1}}{\poqq{q^2}{2n+1}}q^{3 n^2+2n} =\frac {\poqq{-q}{\infty}}{\poqq{q^2}{\infty}}\bigg((-q^{11},- q^{19},q^{30};q^{30})_\infty-q^3(-q,-q^{29},q^{30};q^{30})_\infty\bigg).\label{yyyyyy-newg}
\end{align} 
\end{tl}
\pf It suffices to replace $q$ with $q^2$  in \eqref{yyyyyyg} first.  Then
 \begin{align}
\sum_{n=1}^{\infty}q^{2n^2-2n}\frac{\poqq{x,y}{n}}{\poqq{q^2}{2n-1}}
\left(\frac{1}{xy}\right)^n=(T_1-T_2)\frac {\poqq{q^2/x,q^2/y}{\infty}}{\poqq{q^2,q^2/xy}{\infty}}. \label{yyyyyy-twog}
 \end{align}
 And then set $y\to \infty$ alone. Then \eqref{yyyyyy-twog} reduces to
  \begin{align}
\sum_{n=1}^{\infty}q^{3 n^2-3n}\frac{\poqq{x}{n}}{\poqq{q^2}{2n-1}}\left(-\frac{1}{x}\right)^n
=(T_1-T_2)\frac {\poqq{q^2/x}{\infty}}{\poqq{q^2}{\infty}}, \label{temp-1g}
 \end{align}
 In this form, we further set $x=-q$ and evaluate each $T_i$  in \eqref{temp-1g} via Jacobi's triple product identity  \eqref{triple}. We obtain
 \begin{align*}
&T_1=\sum_{n=-\infty}^\infty\,q^{15 n^2+4n-1}=q^{-1}(-q^{11},- q^{19},q^{30};q^{30})_\infty,\\
 &T_2=\sum_{n=-\infty}^\infty\,q^{15 n^2+14n+2}=q^{2}(-q,- q^{29},q^{30};q^{30})_\infty.
 \end{align*} 
A direct substitution of these results into \eqref{temp-1g} leads us to \eqref{yyyyyy-newg}.
  \qed

From \eqref{yyyyyyg}   we   recover a $q$-identity of Rogers-Ramanujan type modulo $24$.
\begin{tl}[{\rm cf. \cite[(86)]{Slater-1}: mod 24 identity}] \begin{align}
\poq{q}{\infty}\sum_{n=0}^{\infty}\frac{q^{2n^2+2n}}{\poq{q}{2n+1}}
=(-q^{7},-q^{17},q^{24};q^{24})_\infty-q^2(-q,-q^{23},q^{24};q^{24})_\infty.\label{ooo-1g}
\end{align}
\end{tl}
\pf  It only need to take the limitation $x,y\to \infty$ of  \eqref{yyyyyyg} simultaneously. We obtain 
\begin{align*}\poq{q}{\infty}
\sum_{n=1}^{\infty}\frac{q^{2n^2-2n}}{\poq{q}{2n-1}}
=T_1-T_2, 
 \end{align*}
where, by Jacobi's triple product identity \eqref{triple},  
 \begin{align*}
T_1&=\sum_{n=-\infty}^\infty\,q^{12n^2+5n}=(-q^{7},-q^{17},q^{24};q^{24})_\infty,\\
T_2&=\sum_{n=-\infty}^\infty\,q^{12n^2+13n+3}=q^2(-q,-q^{23},q^{24};q^{24})_\infty.
 \end{align*} 
This gives the complete proof of \eqref{ooo-1g}.
\qed

Naturally, we now substitute the Bailey pair $(\alpha_n^{(1)}(1,q),\beta_n^{(1)}(1,q))$ to \eqref{AAA3}, we  recover  Eq.  (3.15.3) in \cite{ramanujan}.
\begin{tl}[Mod 15 identity]
\begin{align}
\sum_{n=0}^{\infty}q^{3n^2/2+3n/2} \frac{(-q;q)_n}{\poq{q}{2n+1}}
\label{220g}=\frac{(-q^4,-q^{11},q^{15};q^{15})_{\infty}-q(-q,-q^{14},q^{15};q^{15})_{\infty}}{\varphi(-q)}.
\end{align}
\end{tl}
\pf It suffices to insert $(\alpha_n^{(1)}(1,q),\beta_n^{(1)}(1,q))$ into \eqref{AAA3}, yielding
\begin{align*}
\poq{q}{\infty}&\sum_{n=1}^{\infty}q^{3n^2/2-3n/2} \frac{(-1;q)_n}{\poq{q}{2n-1}}=\frac{2}{\varphi(-q)} \sum_{n=0}^{\infty} q^{n(n+1) / 2}\tau(n)(1-q^{n})H(n)\\
&=\frac{2\poq{q}{\infty}}{\varphi(-q)}\sum_{n=0}^{\infty}q^{(3n+1)(3n+2) / 2}\tau(3n+1)(1-q^{3n+1})(-1)^{n+1}q^{-3n^2/2-5n/2-1}\\
&+\frac{2\poq{q}{\infty}}{\varphi(-q)}\sum_{n=0}^{\infty} q^{(3n+2)(3n+3) / 2}\tau(3n+2)(1-q^{3n+2})(-1)^{n+1}q^{-3n^2/2-7n/2-2}.
\end{align*}
After a bit simplification, we arrive at
\begin{align*}
\sum_{n=1}^{\infty}q^{3n^2/2-3n/2} \frac{(-1;q)_n}{\poq{q}{2n-1}}&=\frac{2}{\varphi(-q)}\sum_{n=-\infty}^{\infty}q^{15n^2/2-7n/2}-\frac{2q}{\varphi(-q)}\sum_{n=-\infty}^{\infty} q^{15n^2/2-13n/2}\\
&=2\frac{(-q^4,-q^{11},q^{15};q^{15})_{\infty}-q(-q,-q^{14},q^{15};q^{15})_{\infty}}{\varphi(-q)}.\end{align*}
This is, after a bit simplification, in agreement with \eqref{220g}.
\qed

We end our discussions with  the following new summation formula for $q$-series in the sense that it is  different from \eqref{qid18-18}.
 \begin{dl}For any integer $m\geq 0$, it holds
  \begin{align}\sum_{n=1}^m\tau^2(n)\tau(n-m) q^{-n}\frac{\left(q^{m-n+1}; q\right)_{2n}}{\poq{q}{2n-1}}=\left\{\begin{array}{ll}
(q^{2L}-1)q^{3L^2-L-1},&\mbox{if}~~m=3L;\\
(1-q^{4L+2})q^{3L^2+L-1},&\mbox{if}~~m=3L+1;\\
(q^{2 L+2}-1)q^{3L^2+5L+1},&\mbox{if}~~m=3L+2.
\end{array}\right.\label{wangyee-1gg}
\end{align}
 \end{dl}
 \pf It suffices to set $a=1$ in \eqref{wangyee} and substitute $(\alpha_n^{(1)}(1,q),\beta_n^{(1)}(1,q))$ into the resulting identity. Then we have
 \begin{align}
\tau(m) q^{m}\sum_{n=1}^m\left(q^{-m}, q^{1+m} ; q\right)_n \frac{\tau^2(n)q^{-n}}{\poq{q}{2n-1}}=F(m),\label{wangyee-2g}
\end{align}
where, for brief, we write
\[F(m):=\sum_{k=0}^m \tau(k)(1-q^{2k})\frac{H(k)}{\poq{q}{\infty}}.\]
Using the basic relations \eqref{basictworelas} we further simplify
\begin{align*}
\mbox{LHS of \eqref{wangyee-2g}}=
\sum_{n=1}^m\tau^2(n)\tau(n-m) q^{-n}\frac{\left(q^{m-n+1}; q\right)_{2n}}{\poq{q}{2n-1}}.
\end{align*}
On the other hand,  by  \eqref{vvvvvvg}  we can  evaluate $F(m)$ on the right-hand side of  \eqref{wangyee-2g} according as the cases $m\pmod 3$. For instance, let $m=3L$, we have 
\begin{align*}
F(3L)&=\sum_{k=0}^{L-1}\tau(3k+1)(1-q^{6k+2})\frac{H(3k+1)}{\poq{q}{\infty}}+\sum_{k=0}^{L-1}\tau(3k+2)(1-q^{6k+4})\frac{H(3k+2)}{\poq{q}{\infty}}\\
&=\sum_{k=0}^{L-1}q^{3 k^2 - k -1}(1-q^{6k+2})-\sum_{k=0}^{L-1}q^{3 k^2+k-1}(1-q^{6k+4})\\
&=\bigg(\sum_{k=0}^{L-1}q^{ 3 k^2-k-1 }-\sum_{k=0}^{L-1}q^{ 3 k^2+5 k +1 }\bigg)-\bigg(\sum_{k=0}^{L-1}q^{3 k^2+k-1}-\sum_{k=0}^{L-1}q^{3 k^2+7k+3}\bigg)=(q^{2L}-1)q^{3L^2-L-1}.
\end{align*}
For $m=3L+1$  and $m=3L+2$, it is easy to check that 
\begin{align*}
F(3L+1)&=F(3L)+q^{-1 - L + 3 L^2}(1-q^{6L+2})=(1-q^{4L+2})q^{3L^2+L-1},\\
F(3L+2)&=F(3(L+1))=(q^{2L+2}-1)q^{3L^2+5L+1}.
\end{align*} Thus \eqref{wangyee-1gg} is proved.
 \qed
 
\section*{Acknowledgements}
The authors are indebted to  S. O. Warnaar for his helps and insightful comments.  
This  work was supported by the National Natural Science Foundation of
China [Grant No. 12471315].

\end{document}